\documentclass[11pt,letterpaper]{amsart}
\usepackage{amsmath,amssymb,amsthm}
\usepackage{enumerate}
\usepackage{multirow}

\newfont{\Bdd}{msbm10 scaled\magstep1}
\newfont{\footnotesizeBdd}{msbm8 scaled\magstep1}

\newtheorem{myth}{Theorem}

\newtheorem{myex}{Example}

\newtheorem{myre}{Remark}

\def\Gal{\mathop{\rm Gal}\nolimits}

\title{ Thue Equations and CM-Fields}

\author{Yves Aubry}
\address{Institut de Math\'ematiques de Toulon, Universit\'e de Toulon, France and Insitut de Math\'ematiques de Marseille, Aix-Marseille Universit\'e, France}
\email{yves.aubry@univ-tln.fr}

\author{Dimitrios Poulakis}
\address{Aristotle University of Thessaloniki, Department of Mathematics, Thessaloniki 54124, Greece}
\email{poulakis@math.auth.gr}

\begin{document}

\begin{abstract} 
We obtain a polynomial type upper bound for the size of the integral solutions of
  Thue equations $F(X,Y) = b$ defined over a totally real number field $K$,
assuming  that  $F(X,1)$ has  a  root  $\alpha$ such that
  $K(\alpha)$ is a CM-field.
Furthermore, we give an algorithm for the computation of the integral solutions of  such an equation.
 
\end{abstract}

 \maketitle

\section{Introduction}

Let $F(X,Y)$ be an irreducible binary form in $\mbox{\Bdd Z}[X,Y]$ with $\deg F \geq 3$ and $b\in\mbox{\Bdd Z} \setminus\{0\}$.
In  1909,  Thue  \cite{Thue}  proved  that  the  equation
$F(X,  Y)  =  b$ has  only  finitely  many    solutions  $(x,  y)\in\mbox{\Bdd Z}^2$.
Thue's  proof  was ineffective and therefore  does not provide  a method to  determine  the integer
solutions of this equation. Other  non effective  proofs  of  Thue's result    can  be  found  in  \cite[Chap.  X]{Dickson}  and  \cite[Chap. 23]{Mordell}.

In  1968, Baker \cite{Baker}, using his results on linear forms in logarithms
of algebraic numbers, computed an explicit upper bound for the size of the integer
 solutions of Thue equations. Baker's result were improved by several authors
 (see for instance \cite{Bugeaud}, \cite{Gyory1}, \cite{Poulakis1})
  but the bounds  remain  of exponential  type and thus, are not useful to compute integer solutions of such equations.
Nevertheless, computation techniques for the resolution of Thue equations have been developed based on the above results 
\cite{Bilu}, \cite{Hanrot}, \cite{Petho}, \cite{Tzanakis} and the solutions of certain parameterized families
 of Thue equations have been obtained \cite{Heuberger}.
Furthermore, upper bounds for the number of integral solutions of Thue equations have been given \cite{Brinza}, \cite{Evertse}, \cite{B_P_vdP_W}.

In the case where all roots of the polynomial $F(X,1)$ are non real, we have a
polynomial type bound provided by other methods \cite[Theorem 2, page 186]{Mordell}, \cite{Gyory} \cite{Poulakis2}.
Gy\H{o}ry's improvement in \cite[Th\'eor\`eme 1]{Gyory} holds in the case where the splitting field of $F(X,1)$ is a CM-field i.e., 
is a totally imaginary quadratic extension of a totally real number field.
In the same paper, Gy\H{o}ry studied Thue equations defined over a CM-field $L$ and also gave  (\cite[Th\'eor\`{e}me 2]{Gyory})
 a polynomial upper bound for the size of their real algebraic integers solutions in $L$.

In this paper, we consider Thue equations $F(X,  Y)  =  b$ defined over a totally real number field $K$. Simplifying
  Gy\H{o}ry's approach, we obtain (Theorem \ref{thuetheorem}) polynomial type bounds for the size of their integral solutions over $K$,
assuming  that $F(X,1)$ has  a  root  $\alpha$ such that  the  field $K(\alpha)$ is a CM-field.
In case where the splitting field is a CM-field we are in
the situation of \cite[Th\'eor\`eme 2]{Gyory}.
Whenever all  roots of the polynomial $F(X,1)$ are non real and $K\neq \mbox{\Bdd Q}$, we obtain much better bounds than those
 already known \cite{Poulakis2}. Moreover, whenever $F(X,1)$ has a real and a non real root, we obtain  polynomial type bounds  that the Baker's method was
  not able to provide other than exponential bounds. Furthermore, the method of the proof of Theorem 1 provides us with an algorithm 
 for the determination of the solutions of such equations.

We illustrate our result by giving two examples of infinite families of Thue equations $F(X,  Y)  =  b$ satisfying the hypothesis of
 Theorem \ref{thuetheorem}. First we consider Thue equations over some totally real subfields $K$ of cyclotomic fields  $N$
 such that the splitting field $L$ of $F(X,1)$ over $K$ is contained in $N$. In this case, $L$ is an abelian extension
 of $K$.  Next, we give a family of equations $F(X,Y)=b$ such that $F(X,1)$ has  a  root  $\alpha$ for which   $K(\alpha)$ is 
a biquadratic  CM-field. These families contain equations such that $F(X,1)$ has also real roots and so the only method for having
upper bound for the size of theirs solutions is Baker's method which provide only bounds of exponential type.
Finally, we give two examples of determination of solutions of equations satisfying the hypothesis of Theorem 1, by using our algorithm.

\section{New bounds}

We introduce a few notations.
Let $K$ be a number field. We consider the set of absolute values of $K$ by extending  the ordinary absolute value $|\cdot|$ of $\mbox{\Bdd Q}$
and, for every prime $p$, by extending the $p$-adic absolute  value $|\cdot|_{p}$ with $|p|_{p} = p^{-1}$.
Let $M(K)$ be an indexing set of symbols $v$ such that  $|\cdot|_v$, $v\in M(K)$, are all of the
above absolute values of $K$.
 Given such an absolute value $|\cdot|_{v}$ on $K$, we denote by
 $d_{v}$  its local degree. Let ${\bf x} = (x_0:\ldots:x_{n})$ be  a point of the projective space
$\mbox{\Bdd P}^n(K)$ over $K$.  We define the field height $H_{K}(\bf x)$ of $\bf x$ by
$$ H_K({\bf x}) = \prod_{v \in M(K)} \max\{|x_{0}|_{v},\ldots,|x_{n}|_{v}\}^{d_v}.$$
Let $d$ be the degree of $K$. We define  the absolute height $H(\bf x)$ by
$H({\bf x}) = H_K({\bf x})^{1/d}$.
 For $x \in K$, we put $H_K(x) = H_K((1:x))$ and $H(x) = H((1:x))$. If $G \in K[X_1,\ldots,X_m]$,
 then we define the field height $H_{K}(G)$ and the absolute height $H(G)$ of $G$ as
the field height and the absolute height respectively of the point whose coordinates are the coefficients of
$G$ (in any order). For an account of the properties of heights see \cite{Hindry,Lang,Silverman}.
 Furthermore, we denote by $O_K$ and $N_K$ the ring of integers
of $K$ and the norm relative to the extension $K/{\mbox{\Bdd Q}}$, respectively.
Finally, for every $z\in {\mbox{\Bdd C}}$ we denote, as usually, by $\bar{z}$ its complex conjugate.

We prove the following theorem:

 \begin{myth}\label{thuetheorem}
 Let $K$ be  a totally real number field of degree $d$.
 Let $b\in O_K\setminus \{0\}$ and  $F(X,Y)\in O_K[X,Y]$ be a form of degree $n \geq 2$.
Suppose that $F(X,1)$ has   a  root $\alpha$ such that  $K(\alpha)$ is a CM-field.
Then the solutions $(x,y)\in O_K^2$ of $F(X,Y)=b$ satisfy
$$H(x) <  \Omega_1 \ \ {\rm and}\ \  H(y) < \Omega_2$$
for the following values of $\Omega_1$ and $\Omega_2$.
If  the coefficients of $X^n$ and $Y^n$ are $\pm 1$, then
$$\Omega_1 = \Omega_2 =   32  H(b)^{1/n} H(F)^{1+1/n} N_{K}(b)^{2/d}.$$
If only the coefficient of $X^n$ is $\pm 1$, then
$$\Omega_1 = 2^9   H(b)^{1/n} H(F)^{2+1/n} N_{K}(b)^{4/d} \ \  {\rm and} \ \  \Omega_2 = 32  H(b)^{1/n} H(F)^{1+1/n} N_{K}(b)^{2/d}.$$
If both the coefficients of $X^n$ and $Y^n$ are $ \neq \pm 1$, then 
$$\Omega_1 =   2^9   H(b)^{1/n} H(\Gamma)^{2n+1} N_{K}(b)^{4/d} H(a_0)^{n-1} N_K(a_0)^{2(n-1)/d}$$ 
and
$$\Omega_2 =  32  H(b)^{1/n} H(\Gamma)^{n+1} N_{K}(b)^{2/d} H(a_0)^{n-1} N_K(a_0)^{2(n-1)/d},$$
where $a_0$ is the coefficient of $X^n$ and $\Gamma$  a point of the projective space with  $1$ and the coefficients of $F(X,Y)$ as coordinates. 
Furthermore, the  number of integral solutions to the equation $F(X,Y) = b$ is at most 
$$ 72\cdot 4^{dn} N_K(b)^{2n}.$$
In case where  $b$ is a unit of $O_K$,  this number is at most $2wn$, where $w$ is  the number of the roots of unity in $K(\alpha)$. 
\end{myth}

The proof of this result is relied on the following property of CM-fields. A non real algebraic number field $L$ is a CM-field if and only if $L$
is closed under the operation of complex conjugation and complex conjugation
commutes with all the $\mbox{\Bdd Q}$-monomorphisms of $L$ into $\mbox{\Bdd C}$ (\cite{Blanksby}, \cite[Th\'{e}or\`{e}me 1]{Gyory2}, \cite[Lemma 2]{Lou_Oka_Oli}).

When $K =\mbox{\Bdd Q}$ and the splitting field of $F(X,1)$ over $\mbox{\Bdd Q}$ is an abelian totally imaginary extension,  the  hypothesis on complex 
conjugation is obviously satisfied. If  the coefficient of $X^n$ is $\pm 1$, it is interesting to notice  that our bounds are essentially independent
 of the degree of the  form  $F(X,Y)$.
Thus, in case where $H(F)$ and $H(b)$ are not too large, an exhaustive search can provide the  integer solutions we are looking for.

Finally, it should be noticed that in case where $F(X,Y)$ is irreducible we are  in
the situation of \cite{Gyory}. For $K = \mbox{\Bdd Q}$,  \cite[Th\'eor\`eme 1]{Gyory} provides a better upper bound than Theorem 1 and for 
 $K \neq \mbox{\Bdd Q}$,  \cite[Th\'eor\`eme 2]{Gyory} gives similar 
 upper bounds to  Theorem 1. Suppose $F(X,Y)$ is reducible. If   $K = \mbox{\Bdd Q}$ and $F_1(X,Y)$ is a non trivial irreducible factor of $F(X,Y)$
 over $\mbox{\Bdd Z}$ such that the splitting field of $F_1(X,1)$ is of CM-type, then  each integer solution 
$(x,y)$   of  $F(X,Y) = b$ satisfies $F_1(x,y) = b_1$, for some divisor $b_1$ of  $b$. Thus, \cite[Th\'eor\`eme 1]{Gyory} 
applies to this equation and gives better 
explicit bound for $|x|$ and $|y|$ than Theorem 1. Finally, let  $K \neq \mbox{\Bdd Q}$ and   $F_1(X,Y)$ be 
 a non trivial irreducible factor  of $F(X,Y)$ over $O_K$ of degree $\nu$ such that the splitting field of $F_1(X,1)$ is of CM-type. Then each solution 
$(x,y)\in O_K^2$   of  $F(X,Y) = b$ satisfies $F_1(x,y) = b_1$, for some divisor $b_1$ of  $b$. Note that  we do not know the height of $b_1$.
For this we use  \cite[Lemma 3]{Gyory1} which yields a unit $\epsilon \in O_k$ having 
$$H(b_1\epsilon^{\nu}) \leq N_K(b_1)^{1/d} \exp\{c\nu R_K\},$$
 where $c$ is an explicit constant and $R_K$ the regulator of $K$. Thus, we have
 $F_1(\epsilon x,\epsilon y) = b_1\epsilon^{\nu}$ and so, using \cite[Th\'eor\`eme 2]{Gyory}, we obtain an upper bound for $H(x)$ and $H(y)$
 with an  extra factor which is exponential in terms $R_K$  and hence it  is clearly worse than that of Theorem 1.

\section{Examples}

In this section we give two examples in order to illustrate  our result. We denote by $F^*(X,Y)$ the homogenization of a polynomial
$F(X)\in {\mbox{\Bdd C}}[X]$.

\begin{myex}{\rm
Let $p$ be a prime with $p \equiv 1 \ (\bmod\,4)$ and $\zeta_p$ a  $p$-th
primitive root of unity in $\mbox{\Bdd C}$. Then the quadratic
field ${\mbox{\Bdd Q}}(\sqrt{p})$ is a subfield of $\mbox{\Bdd Q}(\zeta_p)$. The field $\mbox{\Bdd Q}(\zeta_p)$ is a cyclic extension of $\mbox{\Bdd Q}$
with Galois group  $\Gal(\mbox{\Bdd Q}(\zeta_p)/\mbox{\Bdd Q}) \simeq (\mbox{\Bdd Z}/p\mbox{\Bdd Z})^*$.

Let $\alpha\in\mbox{\Bdd Z}[\zeta_p]$ be a primitive element of the extension $\mbox{\Bdd Q}(\zeta)/\mbox{\Bdd Q}(\sqrt{p})$ and
  $\alpha_1,\ldots, \alpha_m$, with
$m = (p-1)/2$,  all the distinct conjugates of $\alpha$ over $\mbox{\Bdd Q}(\sqrt{p})$.
The largest real field contained in $\mbox{\Bdd Q}(\zeta_p)$ is $K_p =\mbox{\Bdd Q}(\zeta_p+\bar{\zeta}_p)$ which is a totally real number field.
Let  $\beta\in K_p$ be a primitive element of the extension $K_p/\mbox{\Bdd Q}(\sqrt{p})$ and $\beta_1,\ldots,\beta_n$, where $n = (p-1)/4$,
all the distinct conjugates of $\beta$ over $\mbox{\Bdd Q}(\sqrt{p})$. Then the polynomial
$$F(X) = (X-\alpha_1)\cdots (X-\alpha_m)(X-\beta_1)\cdots (X-\beta_n)$$
belongs to $\mbox{\Bdd Q}(\sqrt{p})[X]$ and has real and non real roots.  Furthermore, we have $\mbox{\Bdd Q}(\sqrt{p})(a_i) = \mbox{\Bdd Q}(\zeta_p)$
which is a CM-field. Consequently, for every non zero
$b\in\mbox{\Bdd Z}[(1+\sqrt{p})/2]$, the Thue equation $F^*(X,Y)=b$ satisfies the hypothesis of  Theorem \ref{thuetheorem}. Note that this equation
satisfies also the hypothesis of \cite[Th\'eor\`eme 2]{Gyory}.

Then, using \cite[Theorem 5.9, page 211]{Silverman} and \cite[Lemma 5.10, page 213]{Silverman},  Theorem \ref{thuetheorem} gives
 the following upper bound
 for the heights of  solutions  $x,y\in\mbox{\Bdd Z}[(1+\sqrt{p})/2]$:
 $$H(x) < 2^{(3p+17)/4} ( H(\alpha)^2H(\beta))^{(3p+1)/6}    H(b)^{4/3(p-1)} N_{\mbox{\footnotesizeBdd Q}(\sqrt{p})}(b)$$
 and
 $$H(y) < 2^{(3p+13)/2}  ( H(\alpha)^2H(\beta))^{5(p-1)/6} H(b)^{4/3(p-1)} N_{\mbox{\footnotesizeBdd Q}(\sqrt{p})}(b)^2.$$
If we consider the particular case where $\Phi_p(X)$ is the $p$-th cyclotomic polynomial, then  \cite[Section 2]{Gyory}
 implies that the maximum of the absolute  heights of
all algebraic integers $x,y \in K_p$ with  $\Phi_p^*(x,y) = 1$  is
$< 2^{(p-1)/2}$.  Theorem  1 improves this result by  yielding the bound  $32$. 
}
\end{myex}

\begin{myex}
{ \rm  Let $d$ be a positive integer $\geq 2$ and $r=m+n\sqrt{d}$, where $m$, $n$ are  integers such that $m>0$ and  $m^2-n^2d>0$.
The minimal polynomial of $ r$ over $\mbox{\Bdd Q}$ is
$$M(X) = X^2-2mX+m^2-dn^2.$$
Then,  the polynomial
$$P(X)= M(-X^2)= X^4+2mX^2+m^2-dn^2$$
is the minimal polynomial of $\sqrt{-r}$  over $\mbox{\Bdd Q}$. Since $m>0$ and $m^2-n^2d>0$, their roots are not real and so, 
$\mbox{\Bdd Q}(\sqrt{d},\sqrt{-r})$ is a CM-field. If $Q(X) \in \mbox{\Bdd Z}[\sqrt{d}][X]\setminus \mbox{\Bdd Z}$, then we put 
$$F(X) = (X^2+(m+n\sqrt{d})Y^2) Q(X).$$
Then for every nonzero $b \in \mbox{\Bdd Z}[\sqrt{d}]$, the Thue equation $F^*(X,Y) = b$  over $K = \mbox{\Bdd Q}(\sqrt{d})$
 satisfies the hypothesis of Theorem 1.  Suppose that $Q(X)$ is monic and $\deg Q = q >0$. By \cite[Remark B.7.4]{Hindry}, we have
 $$H(F) \leq 4H(m+n\sqrt{d})H(Q).$$
 Thus Theorem 1 yields the following upper bounds for the height of integral solutions of the above equations over $K$:
 $$H(x) <     2^{13+1/q}   H(b)^{1/2q} (H(m+n\sqrt{d})H(Q))^{2+1/q} N_{K}(b)^2,$$
 $$H(y) <     2^{7+1/q}  H(b)^{1/2q} (H(m+n\sqrt{d})H(Q))^{1+1/q} N_K(b).$$
Note that in case where the splitting field of $F(X)$ is not a CM-field, \cite[Th\'eor\`eme 2]{Gyory} cannot be applied. Furthermore,
Baker's method can provide only bounds of exponential type.

}
\end{myex}

\section{Proof of Theorem 1}

  Write
 $$F(X,Y) = a_0 (X-\alpha_1 Y)\cdots (X-\alpha_n Y).$$
First, we consider the case where $a_0 = \pm 1$. If $a_0=-1$,  we replace $F(X,Y)$ by $-F(X,Y)$
and $b$ by $-b$ and  then we may suppose that $a_0 = 1$.
By our hypothesis, there is  $j$ such that $K_j = K(\alpha_j)$ is a CM-field.   

Let $x,y\in O_K$ such that $xy \neq 0$ and $F(x,y) = b$. We set $b_j:=x-\alpha_j y$.
Since $K$ is a totally real number field,  we have $x-\bar{\alpha}_j y = \bar{b}_j$. 
Setting $\bar{b}_j = \rho_j b_j$, we obtain the system
$$x-\alpha_j y = b_j, \ \   x-\bar{\alpha}_j y = \rho_j b_j.$$
Eliminating $b_j$ from the above two equations, we get
$$  x = y \frac{\bar{\alpha}_j-\alpha_j\rho_j}{1-\rho_j}.$$
Set
$$A= \frac{\bar{\alpha}_j-\alpha_j\rho_j}{1-\rho_j}.$$
We have
$$H(A) \leq H(\bar{\alpha}_j-\alpha_j\rho_j) H(1-\rho_j) \leq 4 H(\alpha_j)^2 H(\rho_j)^2.$$
Since $\alpha_j$ is not real,  using \cite{Mignotte}, we deduce
 $H(\alpha_j)< 2 H(F)^{1/2}$.
It follows that
$$H(A) \leq 16 H(F) H(\rho_j)^2.$$
Substituting in the equation $F(x,y) = b$ we deduce that
$$y^n F(A,1) = b,$$
and thus
$$H(y)^n \leq H(F(A,1)) H(b)  \leq (n+1)H(F)H(A)^n H(b).$$
Using the bound for $H(A)$ we obtain
\begin{eqnarray}
H(y)^n \leq (n+1)16^n  H(b) H(F)^{n+1} H(\rho_j)^{2n}.
\end{eqnarray}

Next, we shall compute a bound for the height of $\rho_j$. 
We denote by $G_j$ the set of
 $\mbox{\Bdd Q}$-embeddings $\sigma : K_j \rightarrow {\mbox{\Bdd C}}$.
Since $K_j$ is a CM-field, \cite[Th\'{e}or\`{e}me 1]{Gyory2} yields that the complex conjugation
commutes with all the elements of $G_j$. Further,  $K_j$  is closed under the operation of complex conjugation whence we get $\bar{\alpha}_j \in K_j$ and so
$\bar{b}_j\in K_j$. Thus, for every   $\sigma \in G_j$, we have $\sigma(\bar{b}_j) = \overline{\sigma(b_j)}$. It follows that
$$|\sigma(\rho_j) | = \frac{|\sigma(\bar{b}_j)|}{|\sigma(b_j)|} = \frac{|\overline{\sigma(b_j)}|}{|\sigma(b_j)|} = 1.$$
 Let $M_j(X)$ be the minimal polynomial of $\rho_j$ over $\mbox{\Bdd Z}$ and
 $m_j$ its leading coefficient. The elements $\alpha_j$, $\bar{\alpha}_j$ are algebraic integers of $K_j$ and so,
  $b_j$, $\bar{b}_j$ are algebraic integers of $K_j$. It follows that 
 the polynomial
 $$\Pi_j(X) = \prod_{\sigma \in G_j} \sigma(b_j) (X-\sigma(\rho_j))$$
 has integer coefficients. Since $\rho_j$ is a root of $\Pi_j(X)$, we have that  $M_j(X)$ divides $\Pi_j(X)$ and thus we deduce that $m_j$ divides
 $$\prod_{\sigma \in G_j} \sigma(b_j)=N_{K_j}(b_j),$$
 where $N_{K_j}$  is the norm relative to the  extension $K_j/{\mbox{\Bdd Q}}$.
It follows that  $m_j$ divides $N_{K_j}(b_j)$.
As we saw above, all the conjugates $\rho_{j1},\ldots,\rho_{j\mu}$ $(\mu \leq dn)$, of $\rho_j$ are of absolute value 1.  By \cite[page 54]{Lang}, we have
$$
 H_{K_j}(\rho_j) =  m_j  \prod_{i=1}^{\mu} \max\{1,|\rho_{ji}|\} = m_j.
$$
Since $m_j$ divides $N_{K_j}(b_j)$, we have that $m_j$ divides $N_{K_j}(b)$. Thus
\begin{eqnarray}
 H_{K_j}(\rho_j)  \leq N_{K_j}(b).
\end{eqnarray}

Combining the inequalities (1) and (2), we get
$$H(y) \leq  32  H(b)^{1/n} H(F)^{1+1/n} N_{K}(b)^{2/d}.$$
We have
$$H(x) \leq H(A) H(y) \leq 16 H(F) H(\rho_j)^2 H(y)$$
whence we obtain
$$H(x) \leq 2^9   H(b)^{1/n} H(F)^{2+1/n} N_{K}(b)^{4/d}.$$

Suppose now that $a_0 \neq \pm 1$. Write $F(X,1) = a_0 X^n+a_1X^{n-1}+\cdots +a_n $.
Then $a_0\alpha_i$ is a root of  $f(X) = X^n+a_1X^{n-1}+a_2a_0 X^{n-2}+\cdots+a_na_0^{n-1}$ and thus
$a_0\alpha_i$ is an algebraic integer. Denote by $F_1(X,Y)$ the homogenization
of $f(X)$. If $(x,y)\in O_K^2$ is a solution to $F(X,Y)=b$, then $(a_0x,y)$ is  a solution to $F_1(X,Y) = ba_0^{n-1}$.
Denote by $\Gamma$ a point in the projective space with  1 and the coefficients of $F$ as coordinates. Then
we have $H(F_1) \leq H(\Gamma)^n$ and finally, we obtain
$$H(y) \leq  32  H(b)^{1/n} H(\Gamma)^{n+1} N_{K}(b)^{2/d} H(a_0)^{n-1} N_K(a_0)^{2(n-1)/d}$$
and
 $$H(x) \leq 2^9   H(b)^{1/n} H(\Gamma)^{2n+1} N_{K}(b)^{4/d} H(a_0)^{n-1} N_K(a_0)^{2(n-1)/d}.$$
 
Now suppose that $b$ is a unit in $O_K$. Then inequality (2) implies that $H(\rho_j)=1$ and so Kronecker's theorem yields that
$\rho_j$ is a root of unity. 
Let $w$ be the number of the roots of unity in $K_j$. Then we have $w$ choices for $A$ (for the roots of unity  $\neq \pm 1$) and, since $y$ is real, the
equation $y^n F(A,1) = b$ gives us at most $2w$ choices for $y$. Considering also the solutions of the equation with $xy=0$,
we deduce that the number of integral solutions to the
equation $F(X,Y) = b$ is at most $2wn$. 
Finally, suppose that $b$ is not a unit in $O_K$. Using \cite[Lemma 8B]{Schmidt}, we obtain that the number of elements   $\rho_j\in K_j$ with 
$H(\rho_j) \leq  N_K(b)^{1/d}$ is bounded by
$$36\cdot 4^{dn} N_K(b)^{2n}$$
and so the result follows.

\section{An Algorithm}

In this section we give an algorithm for the computation of the integral solutions
to $F(X,Y) = b$ based on the proof of Theorem 1.

\ \\
SOLVE-THUE-1\\
{\em Input:} A totally real number field $K$, a form $F(X,Y)\in O_K[X,Y]$ with $F(X,1)$ monic, $b \in O_K\setminus \{0\}$ and $\alpha$
 a root of $F(X,1)$ such that $K(\alpha)$ is a CM-field.\\
{\em Output:} The integral solutions of $F(X,Y) = b$ over $ K$.
 \begin{enumerate}
 \item Compute the set $\Lambda$  of all the elements $\rho\in K(\alpha)\setminus K$
 having  the absolute values of all theirs conjugates equal to 1 and
 $H_{K(\alpha)}(\rho) | N_{K(\alpha)}(b)$. If $b$ is a unit of $O_K$, then the set $\Lambda$
 consist of all the roots of unity of $K(\alpha)$ which does not belong to $K$.
 
 \item Compute the set $\Xi$ of elements $\xi$ of $K$ of the form:
$$\xi = \frac{\bar{\alpha}-\alpha\rho}{1-\rho},$$
where $\rho\in \Lambda$.

\item Compute the set $S$ of elements $y \in O_K$ such that there is $\xi \in \Xi\cup \{(\bar{\alpha}+ \alpha)/2\}$ with
$$y^n F(\xi,1) = b.$$

\item Output the solutions $(x,y)\in O_K^2$ to $F(X,Y) = b$ with $y\in S$ and  the solutions
$(x,y)\in O_K^2$ with $xy =0$.
\end{enumerate}
 
{\it Proof of Correctness.}  Let $(x,y)\in O_K^2$ be a solution to $F(X,Y)=b$ with $xy\neq 0$. We set $x-\alpha y = \beta$
 and $\rho = \bar{\beta}/\beta$. From the proof of Theorem 1 we have $x = y A,$ where 
$$A= \frac{\bar{\alpha}-\alpha\rho}{1-\rho},$$
 and so $y^n F(A,1)=b$. Since $x,y \in K$, we get $A\in K$. Further, $H_{K(\alpha)}(\rho)$ divides $N_{K(\alpha)}(b)$.
 Let $\beta = \beta_1+\beta_2 i$, with $\beta_1,\beta_2 \in \mbox{\Bdd R}$. 
Suppose that $\rho\in K$. Then $\rho\in \mbox{\Bdd R}$ and $\bar{\beta} = \rho \beta = \rho\beta_1+\rho\beta_2 i$, whence we have 
 $\rho\beta_1 = \beta_1 $  and  $\rho\beta_2 = -\beta_2$.  If $\beta_1 \neq 0$, then $\rho = 1$ and so $\bar{\beta} =  \beta$
 which is a contradiction (since $\alpha$ is not real because $K(\alpha)$ is CM). If $\beta_1 = 0$, then $\rho=-1$.
 In this case we have $A=(\bar{\alpha}+ \alpha)/2$.
 Finally, if $b$ is a unit, then  we have that $H(\rho) =1$
and so $\rho$ is a root of unity in $K(\alpha)$. 
 
 \

Note that there are algorithms for the computation of the elements of a number field of bounded height
\cite{Doyle} and for the computation of roots of unity in a number field \cite[Annexe C]{Molin}. 
As far we know there are not implementations for such algorithms. The other computations
can be carry out by a computational system such as MAGMA or MAPLE.

 \begin{myre}
 {\rm  
 By \cite[page 54]{Lang},   the leading coefficient $m$ of the minimal polynomial of $\rho$ is equal to  $H_{K(\alpha)}(\rho)$.
 Thus, $m\rho \in O_K$.
 }
 \end{myre}

Finally, we give two examples of Thue equations that satisfy the hypothesis of Theorem 1 for which we use the previous algorithm to determine all the integral solutions, the first one having a right-hand side a unit but not the second one.

 \begin{myex}  The only solution of the equation  
 $$(X^2+Y^2) (X^2-\sqrt{2}XY+Y^2) = 1$$
 over $\mbox{\Bdd Z}[\sqrt{2}]$ are $(X,Y)= (\pm 1,0),(0,\pm 1)$.
 \end{myex}
 {\it Proof.} The complex number $i$ is a root of $X^2+1$ and  $K=\mbox{\Bdd Q}(\sqrt{2},i)$ is a CM-field. The roots of unity lying in
 $K\setminus \mbox{\Bdd Q}(\sqrt{2})$  are $\pm i$. Next, we compute 
$$\xi_{\pm} = \frac{-i-i(\pm i)}{1-(\pm i)} = \pm 1.$$
Thus we have the equations $y^4 2(2\pm \sqrt{2}) = 1$. If there is $y\in \mbox{\Bdd Z}[\sqrt{2}]$ satisfying
one of these equations, then $2$ is a unit in $\mbox{\Bdd Z}[\sqrt{2}]$, which is a contradiction since its
norm is not equal to $\pm1$. Furthermore, the solutions $(x,y)\in O_K^2$ with $xy =0$ are $(\pm 1,0)$ and $(0,\pm 1)$.

  \begin{myex}  Consider the form 
 $$F(X,Y)= (X^2+(3-2\sqrt{2})Y^2) (X^2-4XY+\sqrt{2}Y^2)\in \mbox{\Bdd Z}[\sqrt{2}][X,Y].$$
 Then the only solutions of the equation  $F(X,Y)= 3\sqrt{2}-4$    over $\mbox{\Bdd Z}[\sqrt{2}]$ are $(X,Y) = (0,\pm 1)$. 
 \end{myex}
 {\it Proof.}
The given equation belongs to the family of equations of Example 2. Thus, we shall use the above algorithm for the determination of their
solutions.  First, we remark that the equation $F(X,0) = 3\sqrt{2}-4$  has no solution over $\mbox{\Bdd Z}[\sqrt{2}]$ and the only solutions
of  $F(0,Y) = 3\sqrt{2}-4$ over $\mbox{\Bdd Z}[\sqrt{2}]$  are $Y = \pm 1$. 

Set $y = i\sqrt{3-2\sqrt{2}}$ and $K = \mbox{\Bdd Q}({y})$. We  have $N_K(6-4\sqrt{2}) = 16$. 
 We shall compute all the elements $\rho\in K\setminus \mbox{\Bdd Q}(\sqrt{2})$ with $H_K(\rho)|16$
 and having all the absolute values of theirs conjugates equal to 1. 

If $H_K(\rho) = 1$, then $\rho$ is a root of unity in $K\setminus \mbox{\Bdd Q}(\sqrt{2})$.
Since there are not other roots of unity in $K$ than $\pm 1$, we consider the case where 
  $H_K(\rho) > 1$.  Let $H_K(\rho) = 2^{\epsilon}$, where $\epsilon = 1,2$. By Remark 1, we have $\rho= \alpha/2^{\epsilon}$, where $\alpha \in O_K$. 
Using MAGMA, we get the following integral base for $K$: 
$$\omega_0 = 1, \ \ \omega_1 = y, \ \  \omega_2 = \frac{1}{2} (y^2 - 1), \ \   \omega_3 = \frac{1}{4}(y^3 + y^2 - y - 1).$$

Since all the conjugates of $\rho$ have absolute value 1, we obtain the two equalities 
 $$((a_0-2a_1)+a_1\sqrt{2})^2+(2-\sqrt{2})((a_2-2a_3)+a_3\sqrt{2})^2 = 2^{2\epsilon},$$
 and
 $$((a_0-2a_1)-a_1\sqrt{2})^2+(2+\sqrt{2})((a_2-2a_3)-a_3\sqrt{2})^2 = 2^{2\epsilon}.$$
 It follows that	
 \begin{eqnarray}
(a_0-2a_2)^2+2a_2^2+2(a_1-3a_3)^2+2a_3^2 =  2^{2\epsilon}
\end{eqnarray}
 and
  \begin{eqnarray}
 2 a_0 a_2-4 a_2^2+8 a_1 a_3-14a_3^2-a_1^2 =0.
 \end{eqnarray}
From (3) and (4) we deduce that   $a_0$, $a_1$, $a_2$ and $a_3$ are even. Furthermore, we have: $4|a_1$.           
 
Suppose that $\epsilon = 1$.   If $a_2$ or $a_3$ is not zero, then the left-hand side of (3) is $ > 4$ which is a contradiction. Hence 
$a_2 = a_3 = 0$. Similarly, we deduce that $a_1 = 0$. Then $a_0 =\pm 2$ and so $\rho \in \mbox{\Bdd Q}$ which is not the case.

Suppose next  that $\epsilon = 2$.  Putting $a_i = a_i^{\prime}$ $(i=0,1,2,3)$ we have 
 \begin{eqnarray}
(a_0^{\prime}- 2a_2^{\prime})^2+2{a^{\prime}_2}^2+2(a_1^{\prime}-3 a_3^{\prime})^2+2 {a^{\prime}_3}^2 =  4.
\end{eqnarray}
If $a_0^{\prime}-2a_2^{\prime} \neq 0$, then (5) implies that $a_1^{\prime} = a_2^{\prime} = a_3^{\prime} =0$ and so $\rho \in \mbox{\Bdd Q}$ which is
a contradiction. Then $a_0^{\prime}= 2a_2^{\prime}$. If $a_3^{\prime} = 0$, then (5) implies that $a_1^{\prime} = \pm 1$ and so $a_1 =\pm 2$. Since $4|a_1$
we obtain a contradiction. Thus $a_3^{\prime} = \pm 1$. If   $a_1^{\prime}-3a_3^{\prime} = 0$, then   $a_1 = \pm 6$ and so 
4 does not divide $a_1$ which is a contradiction.
Finally suppose that $a_2^{\prime} = 0$. It follows that $a_1^{\prime}-3a_3^{\prime} =\pm 1$.  Thus we have
$$(a_0,a_1,a_2,a_3) = (0,8,0,2),(0,-8,0,-2),(0,4,0,-2),(0,-4,0,2).$$
We see that these values do not satisfy (4). Finally, we have $(\bar{y}+y)/2 = 0$ and we see that 
the equation $Y^4F(0,1)= 3\sqrt{2}-4$ has no solution 
in $\mbox{\Bdd Q}(\sqrt{2})$. The result follows.

 \

{\bf Acknowledgements.} 
This work was done during the visit of the second author at the Department of Mathematics of the University
 of Toulon. The second author wants to thanks this Department for its warm hospitality and fruitful collaboration.
 The authors would like also to thank St\'{e}phane Louboutin and Kalman Gy\H{o}ry for fruitful discussions.


\vskip1cm


\begin{thebibliography}{3}


 \bibitem{Bilu} Y. Bilu and G. Hanrot, Solving Thue equations of high degree,
  {\it J. Number Theory} 60 (1996), 373-392.

 \bibitem{Baker} A. Baker, Contribution to the theory of Diophantine equations, I. On representation of integers
   by  binary  forms, {\it  Philos.  Trans.  Roy.  Sot.  London  Ser.}  A  263  (1968),  173-191.


 \bibitem{Blanksby} P. E. Blanksby and J. H. Loxton,  A Note on the Characterization of CM-fields, {\it J. Austral. Math. Soc.}
 (Series A) 26 (1978), 26-30.
 

 \bibitem{B_P_vdP_W}  B. Brindza, \' A. Pint\' er,  A. van der Poorten and M. Waldschmidt,
     On the distribution of solutions of Thue's equations,
     Number theory in progress, Vol. 1 (Zakopane-Koscielisko, 1997), ed. K. Gy\H{o}ry, H. Iwaniec and J. Urbanowicz, 35-46, de Gruyter, Berlin, 1999.

\bibitem{Brinza}
B. Brindza, J.-H.  Evertse and K. Gy\H{o}ry,
Bounds for the solutions of some Diophantine equations in terms of discriminants,
{\it J. Austral. Math. Soc.} Ser. A 51 (1991), no. 1, 8-26.

 \bibitem{Bugeaud} Y. Bugeaud and K. Gy\H{o}ry, Bounds for the solutions of Thue-Mahler equations
and norm form equations, Acta Arithmetica, LXXIV.3 (1996), 273-292.

\bibitem{Dickson} L.  E.  Dickson,  {\it Introduction  to  the  Theory  of  Numbers,}
    Dover,  New  York,  1957.
    
\bibitem{Doyle}  J. R. Doyle and D. Krumm, Computing algebraic numbers of bounded height, 
{\it Mathematics of Computation} (to appear).


\bibitem{Evertse} J.-H. Evertse, The  number  of  solutions  of  decomposable  form
equations, {\it Inventione Mathematicae,}  122,  (1995) 559-601.


\bibitem{Gyory} K. Gy\H{o}ry, Repr\'{e}sentation des nombres entiers par des formes binaires,
{\it  Publicationes Mathematicae Debrecen,} 24 (3-4)
(1977), 363-375.

\bibitem{Gyory1} K. Gy\H{o}ry and K. Yu, Bounds for the solutions of S-unit equations
and decomposable form equations, {\it Acta Arithmetica,} 123.1 (2006), 9-41.

\bibitem{Gyory2} K. Gy\H{o}ry, Sur une classe des corps de nombres alg\'{e}briques et ses applications, {\it Publicationae Mathematicae,}
	22, (1975), 151-175. 

\bibitem{Hanrot} G. Hanrot, Solving Thue equations without the full unit group. Math.
Comp., 69(229) (2000), 395-405.

\bibitem{Heuberger} C. Heuberger,
Parametrized Thue Equations - A survey, Proceedings of the RIMS symposium ``Analytic Number Theory and Surrounding Areas'', Kyoto, Oct 18-22, 2004, 
RIMS K\^{o}ky\^{u}roku vol. 1511, August 2006, 82-91.

\bibitem{Hindry} M. Hindry and J. H. Silverman, {\it Diophantine Geometry, An Introduction.} Springer-Verlag 2000.

\bibitem{Kappe-Warren} L-C. Kappe and B. Warren, An elementary test for the Galois group of a quartic polynomial,
 {\it The Amer. Math. Monthly}, Vol. 96, No. 2 (1989), 133-137.

\bibitem{Lang} S. Lang, {\it Fundamentals of Diophantine Geometry,} Springer-Verlag, New York - Berlin, 1983.

\bibitem{Lou_Oka_Oli} S. Louboutin, R. Okazaki and M. Olivier, The class number one problem for some non-abelian normal CM-fields, 
{\it Trans. Amer. Math. Soc.} 349 (1997), no. 9, 3657-3678.

\bibitem{Mignotte} M. Mignotte, An inequality of the greatest roots of a polynomial, {\it Elem. Math.} 46 (1991), 85-86.

\bibitem{Molin} P. Molin, {\it Integration numerique et calculs de fonctions} L, Th\`{e}se de Doctorat, Universit\'{e} de Bordeaux I, 2010

\bibitem{Mordell} L. J. Mordell, {\it  Diophantine equations,} Pure and Applied Mathematics, Vol. 30, Academic Press, London-New York 1969.

\bibitem{Petho} A. Peth\H{o}, On the resolution of Thue inequalities, {\it J. Symb. Comput.} 4 (1987), 103-109.

\bibitem{Poulakis1} D. Poulakis, Integer points on algebraic curves with exceptional units, {\it J. Austral. Math. Soc.}
 63 (1997), 145-164.

 \bibitem{Poulakis2}  D. Poulakis, Polynomial Bounds for the Solutions of a Class of Diophantine Equations,
 {\it J. Number Theory} 66, No 2, (1997), 271-281.
 
  \bibitem{Schmidt} W. Schmidt, W. M. Schmidt, {\em Diophantine Approximation and Diophantine Equations,}
  Springer-Verlag 1991.

 \bibitem{Silverman} J. H. Silverman, {\it Arithmetic of Elliptic Curves,} Springer Verlag 1986.

 \bibitem{Thue} A. Thue, \"{U}ber Ann\"{a}herungswerte algebraischer Zahlen,
 {\it J.  Reine  Angew.  Mafh.}  135 (1909),  284-305.

 \bibitem{Tzanakis} N. Tzanakis and B. M. M. de Weger, On the practical solution of the
Thue equation. J. Number Theory, 31(2) (1989), 99-132.





\end{thebibliography}
\end{document}